\documentclass[11pt,reqno]{amsart}
\usepackage{amssymb,latexsym}
\usepackage{graphicx}
\usepackage{fancyhdr}

\theoremstyle{plain}
\numberwithin{equation}{section}
\newtheorem{thm}{Theorem}[section]

\newtheorem{lemma}[thm]{Lemma}

\begin{document}

\setcounter{page}{1}

\title[Fibonacci Sequence, Discrete Probability Distributions and Linear Convolution]{Fibonacci Sequence, Recurrence Relations, Discrete Probability Distributions and Linear Convolution}
\author{Arulalan Rajan, R.~Vittal~Rao}
\address{Centre for Electronics Design and Technology\\
                Indian Institute of Science\\
                Bangalore - 560012\\
                India}
\email{perarulalan@gmail.com, \{mrarul,rvrao\}@cedt.iisc.ernet.in}

\author{Ashok~Rao}
\address{Consultant,~165,~11th~Main,~Saraswathipuram,~Mysore,~India}
\email{ashokrao.mys@gmail.com}

\author{H.~S.~Jamadagni}
\address{Centre for Electronics Design and Technology\\
                Indian Institute of Science\\
                Bangalore - 560012\\
                India}
\email{hsjam@cedt.iisc.ernet.in}

\begin{abstract}
The classical Fibonacci sequence is known to exhibit many fascinating properties. In this paper, we explore the Fibonacci sequence and integer sequences generated by second order linear recurrence relations with positive integer coefficients from the point of view of probability distributions that they induce. We obtain the generalizations of some of the known limiting properties of these probability distributions 
    and present certain optimal properties of the classical Fibonacci sequence in this context. In addition, we also look at the self linear convolution of linear recurrence relations with positive integer coefficients. Analysis of self linear convolution is focused towards locating the maximum in the resulting sequence. This analysis, also highlights the influence that the largest positive real root, of the ``characteristic equation'' of the linear recurrence relations with positive integer coefficients, has on the location of the maximum. In particular, when the largest positive real root is 2, the location of the maximum is shown to depend on whether the sequence length is odd or even.
\end{abstract}

\maketitle

\section{Introduction}
\label{intro}
Starting with any sequence $\{f[n]\}_{n \in \mathbb N}$ of positive real numbers, we can generate a probability distribution on $\{1, 2, \ldots, N\}$ for every $N \in \mathbb N$ as follows:
Let $X$ be a random variable taking values $1, 2, \ldots, N$ such that 
\begin{eqnarray*}\nonumber
P(X=n) = \displaystyle{\frac{f[n]}{\displaystyle{\sum_{k=1}^N}f[k]}}
\end{eqnarray*}
In particular, if $f[n]$ are all positive integers, then we can generate such probability distributions using only integers. In a recent work, Neal \cite{neal} has investigated such probability distributions arising out of the Fibonacci sequence and obtained some limiting properties of these distributions as $N \rightarrow \infty$. \\
While exploring integer sequences from the point of view of generating window functions in the context of designing digital filters for signal processing applications \cite{arul}, we found some generalizations of Neal's limiting properties for a class of integer sequences and an optimal property of the Fibonacci sequence in this class. In this paper, we present these results.

\section{Probability Distributions induced by Fibonacci Sequence}
\label{FibonacciSequence}
 The classical Fibonacci sequence \cite{sloane}, \cite{vajda}, is the positive integer sequence $f[n]$ defined by the second order recurrence relation
  \begin{equation}
  \label{FibEqn}
  f[n] = f[n-1] + f[n-2]
    \end{equation}
with initial conditions $f[1] = 1$, $f[2] = 1$.
The ratio of successive terms of the Fibonacci sequence converges to the well known golden ratio $\varphi$, i.e.,
\begin{equation}
\lim_{k \rightarrow \infty} \frac{f[k+1]}{f[k]} = \varphi = \displaystyle{\frac{1 + \sqrt5}{2}}
\end{equation}
For any positive integer $N$, consider the Fibonacci sequence, $f[1], f[2], \ldots, f[N]$, of length $N$. Using the Fibonacci sequence, we can define a discrete probability distribution as follows.
We define the probability mass function as
\begin{center} 
$p_{_{\bar{N}}}[1], p_{_{\bar{N}}}[2], \ldots, p_{_{\bar{N}}}[N]$ 
\end{center}
 \begin{equation}
 \label{ProbDefX}
 p_{_{\bar{N}}}(n) = \frac{f[n]}{\displaystyle{\sum_{k =1}^N f[k]}} \hspace{0.3in} \forall n = 1, 2, \ldots, N
 \end{equation}
Let $\bar{X}$ \footnote{${\bar{(.)}}$: associated with increasing sequence} be a random variable in $\{1,2, \ldots, N\}$ such that 
\begin{equation}
\label{ProbDef}
P(\bar{X} = n) = p_{_{\bar{N}}}(n); \hspace{0.3in} 1 \leq n\leq N.
\end{equation}

In his work, Neal \cite{neal} showed that, $E(\bar{X})$, the mean or expectation of $\bar{X}$, grows linearly as $N$ increases and hence
\begin{equation}
\label{MeanX}
\lim_{N\rightarrow \infty} E(\bar{X}) = \infty
\end{equation}

Similarly, let $\underline{X}$ be a random variable in $\{1,2, \ldots, N\}$, where the probabilities, $p_{_{\underline{N}}}[n]$, correspond to the Fibonacci sequence, of length $N$, in the decreasing order \footnote{${(\underline{.})}$ : associated with decreasing or the flipped sequence}. We define the probability mass function,
 \begin{equation}
 \label{ProbDefY}
P(\underline{X} = n) = \frac{f[N+1-n]}{\displaystyle{\sum_{k =1}^N f[k]}} = p_{_{\underline{N}}}(n) \hspace{0.2in} ; \hspace{0.2in} 1 \leq n\leq N.
 \end{equation}
For such a random variable, $\underline{X}$, Neal \cite{neal} observed the following limiting properties:
\begin{equation}
\label{MeanY}
\lim_{N\rightarrow \infty} E(\underline{X}) = \varphi + 1
\end{equation}
\begin{equation}
\label{VarianceY}
\lim_{N\rightarrow \infty} Var(\underline{X}) = 2\varphi + 1
\end{equation}
where $E(\underline{X})$ denotes the mean or expectation of $\underline{X}$ and $Var(\underline{X})$ denotes the variance of $\underline{X}$.

The aim of this paper, as mentioned earlier, is to show that such limiting properties are characteristic for integer sequences defined by linear recurrence relations , of any order, with positive integer coefficients. For the sake of simplicity, we shall first discuss the second order recurrence relation.

\section{Second Order Recurrence Relations}
Let $a$, $b$ be any two positive integers and consider the integer sequence defined by the second order recurrence relation
\begin{equation}
\label{SecRecEq}
f[n] = a f[n-1] + b f[n-2] \hspace{0.3in}  \textup{where} \hspace{0.2in}  f[1] = 1; f[2]= 1;
\end{equation}
The characteristic equation of this recurrence relation is given by
\begin{equation}
\label{RecCharEqn}
r^2 - ar - b = 0.
\end{equation}
It is well known that the solution to this Eq.(\ref{RecCharEqn}) is given by
\begin{equation}\nonumber
f[n] =  C[R(a,b)^n - R_s(a,b)^n]
\end{equation}
where, 
\begin{equation}\nonumber
R(a,b) = \displaystyle{\frac{a +\sqrt{a^2 + 4b}}{2}}
\end{equation}
is the positive real root.
The other root of the equation is 
\begin{equation}\nonumber
R_s = - \frac{b}{R(a,b)}
\end{equation}
and 
\begin{equation}\nonumber
C = \displaystyle{\frac{1}{\sqrt{a^2 + 4b}}}
\end{equation}
We note that,
\begin{eqnarray}
\label{opt1}R(a,b) &\geq& \varphi \hspace{0.3in} \textup{for all positive integers} \hspace{0.1in} a, b\\
\label{opt2}R(a,b) &=& \varphi   \hspace{0.3in}\textup{when}\hspace{0.1in} a = b = 1 \textup{(classical Fibonacci sequence)}
\end{eqnarray}
We further observe that the ratio of successive terms of this sequence converges to $R(a,b)$, i.e.,
\begin{equation}
\label{GenSecConverge}
\lim_{k \rightarrow \infty} \frac{f[k+1]}{f[k]} = R(a,b)
\end{equation}
As in section \ref{FibonacciSequence}, for any positive integer $N$, consider the above sequence, of length $N$, i.e.,
\begin{center}
$f[1],f[2], \ldots, f[N]$
\end{center} 
and define probability mass functions,
\begin{center}
$p_{_{\bar{N}}}(n)[1],p_{_{\bar{N}}}(n)[2], \ldots, p_{_{\bar{N}}}(n)[N]$
\end{center} 
as
\begin{equation}
\label{ProbDefRec}
 p_{_{\bar{N}}}(n) = \frac{f[n]}{\displaystyle{\sum_{k =1}^N f[k]}} \hspace{0.3in} \forall n = 1, 2, \ldots, N
\end{equation}
Consider a random integer, $\bar{X}$, in $\{1,2, \ldots, N\}$, such that
\begin{equation}
P(\bar{X} = n) = p_{_{\bar{N}}}(n)
\end{equation}
Similarly let $\underline{X}$ be a random variable in $\{1,2, \ldots, N\}$, (with probabilities corresponding to the decreasing sequence) such that,
\begin{equation}
P(\underline{X} = n) = \frac{f[N+1-n]}{\displaystyle{\sum_{k =1}^N f[k]}} = p_{_{\underline{N}}}(n) \hspace{0.3in} \forall n = 1, 2, \ldots, N
\end{equation}
In the rest of the paper, for short, we shall denote $R(a,b)$ by $R$ and $R_S(a,b)$ by $R_S$.
\section{The Main Theorems on the Limiting Properties}

We prove the following limiting properties.
Analogous to eq.(\ref{MeanX}), we have,
\begin{thm}
 \label{th:MeanLinear}
\begin{equation}\nonumber
\label{MeanNX}
\lim_{N \rightarrow \infty} E(\bar{X}) = \infty
\end{equation}
\end{thm}
Even though $E(\bar{X})$ grows as $N$ increases, its variance still converges. In fact, we have,\\ 
 \begin{thm}
 \label{th:VarianceConstantX}
\begin{equation}\nonumber
\label{VarainceNX}
\lim_{N \rightarrow \infty} Var(\bar{X}) = \frac{R}{(R-1)^2}
\end{equation}
 \end{thm}
In the case of $\underline{X} $, both $E(\underline{X})$ and $Var(\underline{X})$ converge as $N \rightarrow \infty$ and the limit of $Var(\underline{X})$ is the same as that of $Var(\bar{X})$. We have, analogous to eq.(\ref{MeanY}) and eq.(\ref{VarianceY}), the following limits.\\
\begin{thm}
 \label{th:MeanConst}
\begin{equation}\nonumber
\label{MeanNY}
\lim_{N \rightarrow \infty} E(\underline{X}) = \frac{R}{(R-1)}
\end{equation}
 \end{thm}
 \begin{thm}
 \label{th:VarianceConstantY}
\begin{equation}\nonumber
\label{VarianceNY}
\lim_{N \rightarrow \infty} Var(\underline{X}) = \frac{R}{(R-1)^2}
\end{equation}
 \end{thm}
We, then, easily observe using these properties that when $a = b = 1$, the limits of eq.(\ref{MeanNX}), eq.(\ref{MeanNY}) and eq.(\ref{VarianceNY}) reduce to the limits, eq.(\ref{MeanX}), eq.(\ref{MeanY}), eq.(\ref{VarianceY}), obtained for the Fibonacci sequence, by Neal \cite{neal}, described in section \ref{FibonacciSequence}. Using eq.(\ref{opt1}) and eq.(\ref{opt2}), we can prove the following optimality properties of the Fibonacci sequence:
 \begin{thm}
 \label{th:MeanGoldenRatio}
\begin{equation}\nonumber
\label{MeanMax}
\max{\left(\lim_{N\rightarrow \infty} E(\underline{X}): a, b \in \mathbb N \right)} = \lim_{N \rightarrow \infty}\left(E(\underline{X}):a=b=1\right) = \varphi + 1
\end{equation}
\end{thm}
\begin{thm}
\label{th:VarianceGoldenRatio}
\begin{eqnarray}
\label{VarianceMax}
\nonumber\max{\left (\lim_{N\rightarrow \infty}Var(\bar{X}): a,b \in \mathbb N\right)}&=& \max{\left(\lim_{N\rightarrow \infty}Var(\underline{X}): a,b \in \mathbb N\right)}\\
\nonumber &=& \lim_{N \rightarrow \infty}\left(Var(\underline{X}):a=b=1\right) = 2\varphi + 1
\end{eqnarray}
 \end{thm}
We now proceed to prove these theorems.

\section{Some notations used in the following sections}
In order to prove the aforesaid theorems, we make use of certain finite sums, presented in \textbf{Appendix. \ref{Prelims}} and \textbf{Appendix. \ref{FiniteSums}}. We list here, some of the notations used in the subsequent sections.
\begin{center}
\begin{tabular}{cccc}
  Increasing sequence & ~ &~& Decreasing Sequence\\
 $\bar{S} = \displaystyle{\sum_{n=1}^{N}}f[n]$& ~ &~& $\underline{S} = \displaystyle{\sum_{n=1}^{N}}f[N+1-n]$ \\
 $\bar{S}_1 = \displaystyle{\sum_{n=1}^{N}}n f[n]$ & ~ &~& $\underline{S}_1 = \displaystyle{\sum_{n=1}^{N}}n f[N+1-n]$ \\
 $\bar{S}_2 = \displaystyle{\sum_{n=1}^{N}}n^2 f[n]$ &~& ~&$\underline{S}_2 = \displaystyle{\sum_{n=1}^{N}}n^2 f[N+1-n]$ 
\end{tabular}
\end{center}

\section{Proofs of Theorems related to distributions induced by increasing sequences}
\label{Increasing}
We look at the mean and variance of the random variable $\bar{X}$, defined as in section \ref{FibonacciSequence}.\\
\textsc{Theorem {\ref{th:MeanLinear}}:}

 \begin{eqnarray*}
\displaystyle{\lim_{N \rightarrow \infty}} E(\bar{X}) = \infty  
 \end{eqnarray*}

\begin{proof}
We have,
\begin{equation}
\label{mean}
E(\bar{X}) = \frac{\bar{S}_1}{\bar{S}}
\end{equation}
Using eq.(\ref{nfn1}) and eq.(\ref{SumIncrease}), we get
\begin{eqnarray}
E(\bar{X})&\sim& \frac{(NR^{N+2}-(N+1)R^{N+1}+R)(R-1)}{(R-1)^2 R^{N+1}}\nonumber\\
&\sim& \frac{NR-N-1}{R-1}\nonumber\\
\label{expx}
\Rightarrow E(\bar{X}) &\sim& N - \frac{1}{R-1}\\
\nonumber E(\bar{X}) &\sim& \mathcal{O}(N)
\end{eqnarray}\qquad
\end{proof}\
Theorem \ref{th:MeanLinear}, now follows from eq.(\ref{expx}).\\
We also observe the following: \\
 If we now define $\bar Y$ = $(N+1)- \bar{X}$, so that $\bar{Y}$ is also a random variable taking value 1, 2, 3, $\ldots$, $N$ and $\bar{X}$ has been in a sense centralized, then by eq.(\ref{expx}),
 \begin{equation*}
 E(\bar{Y}) \sim 1 + \frac{1}{R-1} = \frac{R}{R-1}
 \end{equation*}
Hence 
\begin{equation*}
\displaystyle{\lim_{N \rightarrow \infty}} E(\bar{Y}) = \frac{R}{R-1} 
\end{equation*}
This $\bar{Y}$ corresponds to $\underline{X}$, which we describe in sec. \ref{Decreasing}.\\
\textsc{Theorem {\ref{th:VarianceConstantX}}:}
\begin{eqnarray*}
\lim_{N \rightarrow \infty} Var(\bar{X}) = \frac{R}{(R -1)^2}
\end{eqnarray*}
\begin{proof}
We have,
\begin{equation}
E(\bar{X}^2) = \frac{\bar{S}_2}{\bar{S}}
\end{equation}
Using eq.(\ref{s2}) and eq.(\ref{SumIncrease}), we get
\begin{equation}
\label{ex2}E(\bar{X}^2) \sim \frac{N^2(R-1)^2 - 2N(R-1) + R+1}{(R-1)^2}
\end{equation}
Using eq.(\ref{expx}) and eq.(\ref{ex2}), we get
\begin{eqnarray}
Var(\bar{X})&=& E(\bar{X}^2)-E(\bar{X})^2\nonumber \\
&\sim &\frac{N^2(R-1)^2 - 2N(R-1) + R+1}{(R-1)^2}-\left(N - \frac{1}{R-1}\right)^2\nonumber\\
\Rightarrow \lim_{N \rightarrow \infty}Var(\bar{X})& = & \frac{(R+1)-1}{(R-1)^2}\nonumber
\end{eqnarray}
Thus, we get the limit of the variance as,
\begin{equation}
\label{varx}
\lim_{N \rightarrow \infty}Var(\bar{X}) = \frac{R}{(R-1)^2}
\end{equation}\qquad
\end{proof}
It is interesting to note from eq.(\ref{expx}) and eq.(\ref{varx}) that, even though the mean increases linearly as the length of the sequence, the variance converges, to a function of $R$.

\section{Proofs of Theorems related to distributions induced by decreasing sequences}
\label{Decreasing}
In this section, we look at the probability distribution generated by reversing the sequence values. We recall that the probability of the random variable $\underline{X}$ taking a value $n$ was defined as, 
\begin{equation}
\label{prob}
P(\underline{X} = n) =\frac{f[N+1-n]}{\displaystyle{\sum_{k=1}^{N}f[N+1-k]}} = p_{_{\underline{N}}}(n)
\end{equation}

\textsc{Theorem {\ref{th:MeanConst}}:}
\begin{eqnarray*}
\displaystyle{\lim_{N \rightarrow \infty}}E(\underline{X}) = 1 + \displaystyle{\frac{1}{R-1}} =  \displaystyle{\frac{R}{R-1}}
\end{eqnarray*}
\begin{proof}
We have,
\begin{equation}
E(\underline{X}) = \frac{\underline{S}_1}{\underline{S}}
\end{equation}
From eq.(\ref{nfnd}), we get
\begin{eqnarray}
E(\underline{X}) &=& \frac{(N+1)\bar{S} - \bar{S}_1}{\bar{S}}\nonumber\\
&=& (N+1) - \frac{\bar{S}_1}{\bar{S}}\nonumber\\
&\sim& (N+1) - (N - \frac{1}{R-1}) \hspace{0.1in}\textup{by eq.(\ref{expx}})\nonumber\\
\label{meand}\Rightarrow  \lim_{N\rightarrow \infty}E(\underline{X}) &=& 1 + \displaystyle{\frac{1}{R-1}} = \displaystyle{\frac{R}{R-1}}
\end{eqnarray}\qquad
\end{proof}

\textsc{Theorem {\ref{th:VarianceConstantY}}:}
\begin{eqnarray*}
\displaystyle{\lim_{N \rightarrow \infty}}E(\underline{X}^2) = 1 + \displaystyle{\frac{3R-1}{(R-1)^2}}
\end{eqnarray*}
\begin{proof}
\begin{equation}
Var(\underline{X}) = E(\underline{X}^2) - E(\underline{X})^2
\end{equation}
Using eq.(\ref{n2fnd}) we first obtain $E(\underline{X}^2)$.
\footnotesize
\begin{eqnarray}
E(\underline{X}^2) &=& \frac{\bar{S}_2}{\underline{S}} \nonumber\\
&\sim& \frac{(N+1)^2\bar{S}}{\bar{S}} - \frac{2(N+1)\bar{S}_1}{\bar{S}} + \frac{\bar{S}_2}{\bar{S}}\nonumber\\
&\sim& (N+1)^2 - 2(N+1)E(\bar{X}) +E(\bar{X}^2)\nonumber\\
&\sim& (N+1)^2 - 2(N+1)\left(N-\frac{1}{R-1}\right) + \left(\frac{N^2(R-1)^2 - 2N(R-1)+(R+1)}{(R-1)^2}\right)\nonumber\\
&\sim& \frac{(R-1)^2 +3R -1}{(R-1)^2}\nonumber
\end{eqnarray}
\normalsize
\begin{equation}
\label{ey2}\Rightarrow \lim_{N \rightarrow \infty}E(\underline{X}^2) = 1 + \displaystyle{\frac{3R-1}{(R-1)^2}}
\end{equation} 
From eq.(\ref{ey2}) and eq.(\ref{meand}), we get 
\begin{eqnarray}
Var(\underline{X}) &=& E(\underline{X}^2) - E(\underline{X})^2\nonumber\\
&\rightarrow& 1 + \frac{3R-1}{(R-1)^2} - \frac{R^2}{(R-1)^2}\nonumber\\ 
\label{VarY}\Rightarrow \lim_{N \rightarrow \infty}Var(\underline{X}) &=& \frac{R}{(R-1)^2}
\end{eqnarray} \qquad
\end{proof}
Thus, we find that, the limit of the variance of the probability distribution with probability defined by eq.(\ref{prob}) to be a simple function $R$, the limit of the ratio of successive elements.

\section{Mean and Variance of discrete probability distributions induced by Fibonacci Sequence}
Now, we look at the probability distribution induced by the standard Fibonacci sequence. For such a distribution, we obtain the value of mean and variance,using eq.(\ref{expx}) and eq.(\ref{ey2}) respectively, by using $a = 1 $, $b = 1$, and $R = \varphi$.
\subsection{Distributions induced by Increasing Sequence}
By eq.(\ref{expx}) and eq.(\ref{varx}), we have 
\begin{equation}
\nonumber E(\bar{X})  \sim \mathcal{O}(N)
\end{equation}
and 
\begin{eqnarray}
\lim_{N \rightarrow \infty}Var(\bar{X}) &=& \frac{R}{(R-1)^2} = \frac{\varphi}{(\varphi-1)^2}\nonumber\\
\label{varf}\Rightarrow \lim_{N \rightarrow \infty}Var(\bar{X}) &=& 2\varphi+1
\end{eqnarray}

Neal \cite{neal} looks at the mean and variance of a Fibonacci distribution. It may be pointed out that the fact that the variance converges in this case was not observed by Neal. The above derivations show that, in this case, though the mean increases linearly as the length of the sequence, the variance converges to a simple function of the golden ratio $\varphi$.

\subsection{Distributions induced by Decreasing Sequence}
\label{DPD}
By eq.(\ref{meand}), we have
\begin{equation}
\nonumber \displaystyle{\lim_{N \rightarrow \infty}}E(\underline{X}) =  \displaystyle{\frac{R}{R-1}} = \varphi + 1
 \end{equation}

From eq.(\ref{VarY}), we have
\begin{eqnarray}
\displaystyle{\lim_{N \rightarrow \infty}}Var(\underline{X}) &=& \frac{R}{(R-1)^2} =\frac{\varphi}{(\varphi-1)^2} \nonumber \\
\label{varyd}\Rightarrow \displaystyle{\lim_{N \rightarrow \infty}}Var(\underline{X}) &=& 2\varphi + 1
\end{eqnarray}
Comparing eq.(\ref{varf}) and eq.(\ref{varyd}), we find that, in the case of the classical Fibonacci sequence, the variance of the distribution converges to the same value, irrespective of whether the probability values are increasing or decreasing.

In the following sections, we make a few observations and remarks related to the classical Fibonacci sequence.

\section{Ratio of Variance to Mean}
We observe the curious fact that, from eq.(\ref{MeanNY}) and eq.(\ref{VarianceNY}),
\begin{equation}
\lim_{N\rightarrow \infty} \frac{Var(\underline{X})}{E(\underline{X})} = \frac{1}{R-1} 
\end{equation}
In the case $a = b = 1$, that is, in the case of the distribution induced by the Fibonacci sequence, since $R = \varphi$
\begin{equation}
\lim_{N\rightarrow \infty} \frac{Var(\underline{X})}{E(\underline{X})} = \frac{1}{\varphi-1} = \varphi
\end{equation}
Thus the ratio of the variance of $\underline{X}$ to the expectation of $\underline{X}$ also converges to the Golden Ratio, $\varphi$.

\section{Optimal Probabilistic limit properties of the Fibonacci Sequence}
In this section, we establish certain optimal probabilistic limit properties of the classical Fibonacci sequence. 

\subsection{Mean}
\textsc{Theorem {\ref{th:MeanGoldenRatio}}:}
\begin{eqnarray*}
\label{MeanMax1}
\max{\left(\lim_{N\rightarrow \infty} E(\underline{X}): a, b \in \mathbb N\right)} = \lim_{N \rightarrow \infty}\left(E(\underline{X}) : a,b = 1\right) = \varphi + 1
\end{eqnarray*}
\begin{proof}
 We know that 
\begin{eqnarray}
 R &\geq& \varphi\nonumber\\
 \Rightarrow R-1 &\geq& \varphi -1 \nonumber\\
\label{compare}\Rightarrow \frac{1}{R-1} &\leq& \frac{1}{\varphi-1}
 \end{eqnarray}
 But $\varphi$ = $\displaystyle{\frac{1}{\varphi-1}}$. \\Hence, from eq.(\ref{compare}), we get
\begin{eqnarray}
\frac{1}{R-1} \leq \varphi\nonumber\\
\label{meancomp}\Rightarrow 1 + \frac{1}{R-1} &\leq& 1+\varphi\\
\displaystyle{\frac{R}{R-1}} &\leq& \varphi + 1 \nonumber\\
\Rightarrow \lim_{N\rightarrow \infty}\left(E(\underline{X}):a,b\in \mathbb N\right) &\leq& \lim_{N\rightarrow \infty}\left(E(\underline{X}): a = b = 1\right)
\end{eqnarray}
Thus, \\
$\max{\left(\displaystyle{\lim_{N\rightarrow \infty}} E(\underline{X}): a, b \in \mathbb N\right)} = \displaystyle{\lim_{N \rightarrow \infty}}\left(E(\underline{X}) : a,b = 1\right) = \varphi + 1$
\end{proof}

\subsection{Variance}

Among the family of distributions, induced by second order linear recurrence relations with positive integer coefficients, $a$ and $b$, the Fibonacci sequence gives the maximum limit for $Var(\underline{X})$, i.e.\\
\textsc{Theorem {\ref{th:VarianceGoldenRatio}}:}
\begin{eqnarray}
\label{VarianceMax1}
\max{\left(\lim_{N\rightarrow \infty}Var(\bar{X}): a,b \in \mathbb N\right)}= \max{\left(\lim_{N\rightarrow \infty}Var(\underline{X}): a,b \in \mathbb N\right)}\\
 = \lim_{N \rightarrow \infty}\left(Var(\underline{X}):a,b=1\right) = 2\varphi + 1
\end{eqnarray}
. 
We prove this as follows:
\begin{proof}
\begin{eqnarray}
 R &\ge& \varphi\nonumber\\
 \Rightarrow \frac{1}{R} &\le& \frac{1}{\varphi}\nonumber\\
\Rightarrow -\frac{1}{R} &\ge& - \frac{1}{\varphi}\nonumber\\
\Rightarrow 1-\frac{1}{R} &\ge& 1 - \frac{1}{\varphi}\nonumber\\
\Rightarrow R\left(1-\frac{1}{R}\right)^2 &\ge& \varphi\left(1 - \frac{1}{\varphi}\right)^2\nonumber\\
\label{vareq}\Rightarrow \frac{1}{R\left(1-\frac{1}{R}\right)^2} &\le& \frac{1}{\varphi\left(1 - \frac{1}{\varphi}\right)^2}
\end{eqnarray}
From eq.(\ref{vareq}), we get 
\begin{eqnarray}
\frac{R}{(R-1)^2} \le \frac{\varphi}{(\varphi-1)^2} = 2\varphi + 1\nonumber\\
\label{varmaxi}\Rightarrow \frac{R}{(R-1)^2} &\le& 2\varphi + 1\\
\Rightarrow \label{varmaxi1}\lim_{N\rightarrow \infty} \left(Var(\underline{X}):a,b\in \mathbb N\right) &\leq& \lim_{N\rightarrow \infty}\left(Var(\underline{X}): a = b = 1\right)
\end{eqnarray}
From eq.(\ref{varmaxi1}), we conclude that the limit variance of discrete probability distributions, induced by decreasing sequence and generated using second order linear recurrence, is maximum for the one generated using Fibonacci sequence.
\end{proof}

In the next section, we consider the probability distributions generated by higher order recurrence relation with positive integer coefficients.

\section{Higher order Recurrence relation with positive integer coefficients}
From Section.\ref{Increasing} and Section.\ref{Decreasing}, one can observe that the mean and variance of the probability distributions, induced by general second order linear recurrence relation with positive integer coefficients, depend only on the largest positive real root, $R$, of the associated ``characteristic equation''. It turns out that, even in the case of any other higher order recurrence relation with positive integer coefficients, it is the largest positive real root that determines the mean and variance. We now proceed to prove the same. \\
Let 
\begin{equation}
 f[n] = a_{n-1}f[n-1] + \ldots + a_{n-k}f[n-k] \hspace{0.3in} \forall a_i \in \mathbb{N}
\end{equation}
be a $k^{th}$ order recurrence relation. The characteristic equation of this recurrence relation is given by,
 \begin{equation}
 \label{higherordgenrec}
    f(x) = x^k - a_{n-1}x^{k-1}- \ldots - a_{n-k}
 \end{equation}
\begin{lemma}
\label{LemCharEq}
Eq.(\ref{higherordgenrec}), has atleast one real positive root and if $\alpha$ is the largest positive real root, then all the other roots (real and complex) lie within the disc $|z| < \alpha$.
\end{lemma}

Using Lemma \ref{LemCharEq}, it is easy to infer the following:\\
With $R = \alpha$,
\begin{itemize}
\item The same asymptotic relations as in eq.(\ref{meand})and eq.(\ref{VarY}) hold. 
\item The same limiting properties as in eq.(\ref{MeanY}) and eq.(\ref{VarianceY}) hold.
\item The optimality properties of the Fibonacci sequence as in Theorem \ref{th:MeanGoldenRatio} and in Theorem \ref{th:VarianceGoldenRatio} hold good. The reason is that, the largest positive real root, $\alpha > \varphi$, the golden ratio.
\end{itemize}
Thus, Theorem \ref{th:MeanGoldenRatio} can be generalized as 
\begin{eqnarray*}
\max{\left(\lim_{N\rightarrow \infty} E(\underline{X}): a_i, i \in \mathbb N \right)} = \lim_{N \rightarrow \infty}\left(E(\underline{X}):a_1 = a_2 = 1, a_{i}=0, \forall i> 2\right) = \varphi + 1
\end{eqnarray*}
and Theorem \ref{th:VarianceGoldenRatio} can be further generalized as 
\begin{eqnarray*}
\max{\left(\lim_{N\rightarrow \infty}Var(\bar{X}): a_i, i \in \mathbb N\right)}&=& \max{\left(\lim_{N\rightarrow \infty}Var(\underline{X}): a_i, i \in \mathbb N\right)}\\
\nonumber &=& \lim_{N \rightarrow \infty}\left(Var(\underline{X}):a_1=a_2=1, a_{i}=0,\forall i > 2\right) = 2\varphi + 1
\end{eqnarray*}

Now, we proceed to prove Lemma \ref{LemCharEq}.
\begin{proof}
By Descartes rule of signs, there exists a positive root. Since, \\ 
 \begin{center}
 $ f(0) < 0$, $f(\infty) > 0$ \\
 \end{center}
Let $\alpha$ be the largest positive root. \\

\subsection{$\alpha$ is simple root}
Claim: The root is simple. \\
Let $f(x)$ be as in eq.(\ref{higherordgenrec}). Since $\alpha$ is a root, we have
 \begin{equation}
 \label{alpharoot}
\alpha^k - a_{n-1}\alpha^{k-1} + \ldots + a_{n-k} = 0
 \end{equation}
Suppose, if $\alpha$ is a multiple root, then we must have
\begin{equation}
 f'(\alpha) = 0 
\end{equation}
\begin{equation}
\Rightarrow k\alpha^{k-1}-(k-1)a_{n-1}\alpha^{k-2} + \ldots + a_{n-k+1} = 0 
\end{equation}
Multiplying by $\alpha$, we get
\begin{eqnarray}
 k\alpha^{k} - (k-1)a_{n-1}\alpha^{k-1} + \ldots + a_{n-k+1}\alpha = 0\nonumber\\
 \Rightarrow  k\alpha^{k} &=& (k-1)a_{n-1}\alpha^{k-1} + \ldots + a_{n-k+1}\alpha\nonumber\\
 \Rightarrow   &<& k\left( a_{n-1}\alpha^{k-1} + \ldots +a_{n-k+1}\alpha \right)\nonumber\\
 \Rightarrow     &<& k\left( a_{n-1}\alpha^{k-1} + \ldots +a_{n-k+1}\alpha + a_{n-k}\right)\nonumber\\
 \Rightarrow \alpha^{k} &<&  a_{n-1}\alpha^{k-1} + \ldots + a_{n-k+1}\alpha + a_{n-k}\nonumber\\
 \Rightarrow\alpha^{n} &<& \alpha^{n} \nonumber
\end{eqnarray}
This contradicts eq.(\ref{alpharoot}). Hence $\alpha$ must be a simple root.

\subsection{All roots lie within circle of radius $|\alpha|$}
Now, we prove that all other roots of eq.(\ref{higherordgenrec}) lie within the circle of radius $\alpha$. Consider eq.(\ref{higherordgenrec}), the characteristic equation of the corresponding $k^{th}$ order recurrence relation with positive integer coefficients. Since $\alpha$ is a root, 
   \begin{equation}\nonumber
     f(\alpha) = 0
   \end{equation}
Also, we have
 \begin{eqnarray}
  \displaystyle{\lim_{x\rightarrow \infty}} f(x) &=& \infty \nonumber\\
  \Rightarrow  \label{functionroots} f(x) &\ge& 0, \hspace{0.3in} \forall x \ge \alpha \\
  \textup{and}\label{functiongreaterroots}\hspace{0.1in} f(x) &>& 0  \hspace{0.3in} \forall x > \alpha           
 \end{eqnarray}
Let $\beta$ be any root (real or complex) of eq.(\ref{higherordgenrec}). Clearly, $\beta$ cannot be zero, because,
 \begin{equation}\nonumber
 f(0) = -a_{n-k} < 0
  \end{equation}
Let $b$ = $|\beta|$.\\
Claim: $b$ $\le$ $\alpha$\\
Reason: Suppose $b$ $>$ $\alpha$.\\
 Since $\beta$ is a root, we have
 \begin{eqnarray}
 f(\beta) = 0 \Rightarrow \beta^{k}  &=&  a_{n-1}\beta^{k-1} + \ldots + a_{n-k}\beta^{n-k} \nonumber\\
              \Rightarrow |\beta|^{k}  &=&  |a_{n-1}\beta^{k-1} + \ldots + a_{n-k}\beta^{n-k}|\nonumber\\
             \Rightarrow \beta^{k}  &\leq&  a_{n-1}|\beta|^{k-1} + \ldots + a_{n-k}|\beta|^{n-k}\nonumber\\
              \Rightarrow \beta^{k}-a_{n-1}|\beta|^{k-1} - \ldots - a_{n-k}|\beta|^{n-k} &\leq& 0 \nonumber\\
              \Rightarrow b^{k}-a_{n-1}b^{k-1} - \ldots - a_{n-k}b^{n-k} &\leq& 0 \nonumber\\
  \label{functionlessthanzero} \Rightarrow f(b) &\leq& 0
 \end{eqnarray}
However, since $b$ $>$ $\alpha$, by eq.(\ref{functiongreaterroots}), $f(b)$ $>$ 0. Thus, we get a contradiction and hence $b$ $\le$ $\alpha$.

Claim: $b$ = $\alpha$ \\
Since $\beta$ is a root, we have,
\begin{eqnarray}
\beta^k = a_{n-1}\beta^{k-1} + \ldots + a_{n-k}\beta^{n-k}\nonumber\\
\Rightarrow |\beta|^k &=& |a_{n-1}\beta^{k-1} + \ldots + a_{n-k}\beta^{n-k}|\nonumber
\end{eqnarray}
Since we have assumed that $|\beta|$=$b$=$a$, 
\begin{eqnarray}
\alpha^k = |a_{n-1}\beta^{k-1}+ \ldots + a_{n-k}\beta^{n-k}|\nonumber\\
a_{n-1}\alpha^{k-1}+ \ldots + a_{n-k} &=&  |a_{n-1}\beta^{k-1}+ \ldots + a_{n-k}| \hspace{0.2in}\textup{by \hspace{0.1in} eq.(\ref{alpharoot})} \\
a_{n-1}|\beta|^{k-1}+ \ldots + a_{n-k} &=&  |a_{n-1}\beta^{k-1}+ \ldots + a_{n-k}|
\end{eqnarray}
$\Rightarrow a_n-1\beta^{k-1}, \ldots, a_{n-k}$ are all complex numbers in the same direction. \\
$\Rightarrow a_{n-k}, a_{n-k+1} \beta$ are real positive\\
$\Rightarrow \beta$ is real positive\\
$\Rightarrow \beta = \alpha$.\\
Thus, based on the above claims, we find that 
\begin{itemize}
\item the only root on the circle $|z|$ = $\alpha$ is $z$ = $\alpha$ and there is no root outside this circle
\item all roots not equal to $\alpha$ are inside the circle $|z|$ = $\alpha$.
\end{itemize}

\subsection{$\alpha$ is strictly greater than the golden ratio $\varphi$}
We prove that the largest positive real root of the characteristic equation, corresponding to higher order linear recurrence relation with integer coefficients, is greater than the golden ratio $\varphi$.
\begin{eqnarray}
\alpha^{k} &=&  a_{n-1}\alpha^{k-1} + \ldots + a_{n-k}\alpha^{n-k} \nonumber\\
\textup{Since all terms on the RHS are greater than 0} \nonumber\\
\Rightarrow \alpha^{k} &>&  a_{n-1}\alpha^{k-1} + a_{n-2}\alpha^{k-2}\nonumber \\
\Rightarrow \alpha^{2} &>&  a_{n-1}\alpha + a_{n-2}\nonumber
\end{eqnarray}
As $a_{n-1}$, $a_{n-2}$ are positive integers, we see that,
\begin{eqnarray}
\Rightarrow \alpha^{2} &>&  \alpha + 1\nonumber\\
\label{alphaphicompare}\Rightarrow \alpha^{2} -  \alpha - 1\ &>& 0
\end{eqnarray}
From the fact that, $\varphi$, the golden ratio, is the largest positive root of $x^2-x-1 = 0$ and that $x^2-x-1 > 0$ for $x > \varphi$, it follows that
$\alpha$ is always greater than $\varphi$.
\end{proof}

\section{Convergence in Distribution}
It must be noted that the distribution  $\{p_{_{\underline{N}}}(n)\}, 1 \leq n \leq N$, converges to a geometric distribution on $\mathbb N = \{1,2,3, \ldots\}$. More precisely, we have the following:\\
Let $Y$ be a random integer in $\mathbb N$, such that
\begin{equation}
\label{pn} p(n) = P(Y=n) = \rho^n \left(\displaystyle{\frac{1-\rho}{\rho}}\right)
\end{equation}
where $\rho = \displaystyle{\frac{1}{R}}$.
We have, 
\begin{equation}
\label{sumrho}
\displaystyle{\sum_{n=1}^{\infty}\rho^n} = \displaystyle{\frac{\rho}{1-\rho}}
\end{equation}
and hence
\begin{equation} \nonumber
\displaystyle{\sum_{n=1}^{\infty}n\rho^{n-1}} = \frac{\mathrm{d}}{\mathrm{d}\rho}\left(\displaystyle{\frac{\rho}{1-\rho}}\right) =  \displaystyle{\frac{1}{(1-\rho)^2}}
\end{equation}
Thus,
\begin{equation}\nonumber
\displaystyle{\sum_{n=1}^{\infty}n\rho^{n}} = \displaystyle{\frac{\rho}{(1-\rho)^2}}
\end{equation}
Consequently,
\begin{eqnarray}
E(Y) = \displaystyle{\sum_{n=1}^{\infty}n\rho^{n}}\left(\displaystyle{\frac{1-\rho}{\rho}}\right) &=& \displaystyle{\frac{\rho}{(1-\rho)^2}}\left(\displaystyle{\frac{1-\rho}{\rho}}\right)\nonumber\\
 &=& \displaystyle{\frac{1}{1-\rho}}\nonumber\\
 &=& \displaystyle{\frac{1}{1-\displaystyle{\left(\frac{1}{R}\right)}}}\nonumber\\
\Rightarrow E(Y) &=& \displaystyle{\frac{R}{R-1}}
\end{eqnarray}
Similarly we can see that 
\begin{equation}
Var(Y^2) = \displaystyle{\frac{R}{(R-1)^2}}
\end{equation}
Thus, we see that,
\begin{eqnarray}
\displaystyle{\lim_{N \rightarrow \infty}}E(\underline{X}) = E(Y)\nonumber\\
\displaystyle{\lim_{N \rightarrow \infty}}Var(\underline{X}) = Var(Y)\nonumber
\end{eqnarray}
We now show that the random variable $\underline{X}$ converges to the random variable $Y$ in distribution. \\
We have for every $n \in \mathbb N$,
\begin{eqnarray}
p_{_{\underline{N}}}(n) &\sim& \left(\displaystyle{\frac{R^{N+1-n}}{R^{N+1}- R}}\right)(R-1)\nonumber\\
&=& \left(\displaystyle{\frac{R^{-n}}{1-R^{-N}}}\right)(R-1)\nonumber\\
&=& \left(\displaystyle{\frac{\rho^{n}}{1-\rho^{N}}}\right)\left(\displaystyle{\frac{1}{\rho}}-1\right)\nonumber\\
&=& \left(\displaystyle{\frac{\rho^{n}}{1-\rho^{N}}}\right)\left(\displaystyle{\frac{1-\rho}{\rho}}\right)\nonumber
\end{eqnarray}
Thus,
\begin{equation}\nonumber
p_{_{\underline{N}}}(n) \sim \left(\displaystyle{\frac{\rho^{n}}{1-\rho^{N}}}\right)\left(\displaystyle{\frac{1-\rho}{\rho}}\right)
\end{equation}
Since $\displaystyle{\lim_{N \rightarrow \infty}}\rho^N = 0$ as $\rho = \displaystyle{\frac{1}{R}} < 1$, we have
\begin{equation}
\label{ProbConverge}
\displaystyle{\lim_{N \rightarrow \infty}}p_{_{\underline{N}}}(n) = \rho^{n}\left(\displaystyle{\frac{1-\rho}{\rho}}\right)
\end{equation}
From eq.(\ref{pn}) and eq.(\ref{ProbConverge}) we see that
\begin{equation}
\displaystyle{\lim_{N \rightarrow \infty}}p_{_{\underline{N}}}(n) = p(n) \hspace{0.2in} \forall n \in \mathbb N
\end{equation}
and thus $\underline{X}$ converges to $Y$ in distribution.

\section{Convergence of Moments}
In the previous subsection, we observed that the distributions generated by second order linear recurrence relation converge to a geometric distribution, i.e., 
\begin{itemize}
\item $\underline{X}$ converges to $Y$ in distribution and 
\item $E(\underline{X})$ converges to $E(Y)$  
\item $E(\underline{X}^2)$ converges to $E(Y^2)$ (and hence $Var(\underline{X})$ converges to $Var(Y)$).
\end{itemize}
Now, we prove that all the moments of $\underline{X}$ converge to the corresponding moment of $Y$, as $N \rightarrow \infty$.
\begin{proof}
For any positive integer $k$, we have,
\begin{equation}
\label{Moments}
E(\underline{X}^k)  \sim  \displaystyle{\frac{\displaystyle{\sum_{n=1}^{N}}n^k R^{N+1-n}}{\displaystyle{\sum_{n=1}^{N}}R^{N+1-n}}}
\end{equation}
We have,
\begin{equation}
\label{sumR1}
\displaystyle{\sum_{n=1}^{N}}R^{N+1-n} = \displaystyle{\frac{R^{N+1}-R}{R-1}}
\end{equation}
Also,
\begin{equation}
\label{nsumR1}
\displaystyle{\sum_{n=1}^{N}}n^kR^{N+1-n} = R^{N+1}\displaystyle{\sum_{n=1}^{N}n^kR^{-n}}
\end{equation}
Substituting eq.(\ref{sumR1}) and eq.(\ref{nsumR1}) in eq.(\ref{Moments}), we get
\begin{eqnarray}
E(\underline{X}^k) &\sim& \left(\displaystyle{\frac{R^{N+1}\displaystyle{\sum_{n=1}^{N}n^kR^{-n}}}{R^{N+1}-R}}\right)(R-1)\nonumber\\
&=& \left(\displaystyle{\frac{R^{N+1}\displaystyle{\sum_{n=1}^{N}n^kR^{-n}}}{R^{N+1}(1-R^{-N})}}\right)(R-1)\nonumber\\
&=& \displaystyle{\frac{\displaystyle{\sum_{n=1}^{N}}n^k\rho^{n}}{1-R^{-N}}}\left(\frac{1-\rho}{\rho}\right)\nonumber
\end{eqnarray}
Hence, as $N \rightarrow \infty$,
\begin{eqnarray}
E(\underline{X}^k) &\longrightarrow& \displaystyle{\sum_{n=1}^{N}n^k\rho^{n}}\left(\frac{1-\rho}{\rho}\right) \hspace{0.2in} \textup{because}\hspace{0.2in} R^{-N}\rightarrow 0\nonumber\\
&=& E(Y^k)\nonumber
\end{eqnarray}
Thus, we have for every non-negative integer $k$,
\begin{center}
 $E(\underline{X}^k) \longrightarrow E(Y^k),$ 
\end{center}
showing that all the moments of $\underline{X}$ converge to the corresponding moments of $Y$. 
\end{proof}

\subsection{Limits of Moments in  the case of classical Fibonacci sequence}
The moments of the probability distribution induced by classical Fibonacci sequence converge to a linear function of the Golden ratio, $\varphi$.
\begin{proof}
We have,
\begin{equation}\nonumber
\displaystyle{\sum_{n=1}^{\infty}}x^n = \displaystyle{\frac{x}{1-x}} \hspace{0.3in} \textup{for} \hspace{0.2in} |x| < 1.
\end{equation}
which on differentiation gives
\begin{equation}
\label{ndiff1}
\displaystyle{\sum_{n=1}^{\infty}}nx^{n-1} = \displaystyle{\frac{1}{(1-x)^2}}
\end{equation}
From eq.(\ref{ndiff1}), we get
\begin{equation}
\label{ndiff2}
\displaystyle{\sum_{n=1}^{\infty}}nx^{n} = \displaystyle{\frac{x}{(1-x)^2}}
\end{equation}
Then, 
\begin{equation}
\label{ndiff3}
\displaystyle{\sum_{n=1}^{\infty}}nx^{n}\left(\displaystyle{\frac{1-x}{x}}\right) = \displaystyle{\frac{1}{(1-x)}}
\end{equation}
Similarly, differentiating eq.(\ref{ndiff2}), we get
\begin{eqnarray}
\displaystyle{\sum_{n=1}^{\infty}}n^2x^{n-1} &=& \displaystyle{\frac{(1-x)^2 + 2x(1-x)}{(1-x)^4}}\nonumber\\
\label{ndiff4}= \displaystyle{\frac{(1-x)+2x}{(1-x)^3}}&=& \displaystyle{\frac{1+x}{(1-x)^3}}
\end{eqnarray}
and hence,
\begin{equation}\nonumber
\displaystyle{\sum_{n=1}^{\infty}}n^2x^{n} = \displaystyle{\frac{(1+x)x}{(1-x)^3}}
\end{equation}
and thus,
\begin{equation}
\label{ndiff5}
\displaystyle{\sum_{n=1}^{\infty}}n^2x^n\left(\frac{1-x}{x}\right) = \displaystyle{\frac{1+x}{(1-x)^2}}
\end{equation}
We shall now show, by induction, that in general, for any positive integer $k$,
\begin{equation}
\label{genind}
\displaystyle{\sum_{n=1}^{\infty}}n^k x^n \left(\frac{1-x}{x}\right) = \displaystyle{\frac{P_{_{k-1}}(x)}{(1-x)^k}}
\end{equation}
where $P_{_{k-1}}(x)$ is a polynomial in $x$ of degree $k-1$.
We have from eq.(\ref{ndiff3}) and eq.(\ref{ndiff4}), that the result eq.(\ref{genind}) is true for $k = 1$ and $k = 2$. Thus, it is enough to show that if eq.(\ref{genind}) is true for some $k = r$, then it is also true for $k = r+1$. Suppose eq.(\ref{genind}) is true for $k = r$, then
\begin{equation}\nonumber
\displaystyle{\sum_{n=1}^{\infty}}n^r x^n \left(\frac{1-x}{x}\right) = \displaystyle{\frac{P_{_{r-1}}(x)}{(1-x)^r}}
\end{equation}
Hence, we get,
\begin{equation}
\label{ndiff6}
\displaystyle{\sum_{n=1}^{\infty}}n^r x^n = \frac{xP_{_{r-1}}(x)}{(1-x)^r}
\end{equation}
On differentiating eq.(\ref{ndiff6}), we get,
\begin{eqnarray}
\displaystyle{\sum_{n=1}^{\infty}}n^{r+1}x^{n-1} &=& \frac{(1-x)^{r+1}\{P_{_{r-1}}(x)+P_{_{r-1}}'(x)\}+ xP_{_{r-1}}(x)r(1-x)^r}{(1-x)^{2r+2}}\nonumber\\
&=& \frac{(1-x)\{P_{_{r-1}}(x) + xP_{_{r-1}}'(x)\}+ xP_{_{r-1}}(x)r}{(1-x)^{r+2}}\nonumber\\
\label{ndiff6sol}&=& \frac{P_{_r}(x)}{(1-x)^{r+2}}
\end{eqnarray}
where 
\begin{equation}
\label{ndiff7}
P_{_r}(x) = (1-x)\{P_{_{r-1}}(x) + xP_{_{r-1}}'(x)\} + xP_{_{r-1}}(x)r
\end{equation}
is a polynomial of degree $r$. From eq.(\ref{ndiff6sol}), we get,
\begin{equation}\nonumber
\displaystyle{\sum_{n=1}^{\infty}}n^{r+1}x^n = \frac{xP_{_r}(x)}{(1-x)^{r+2}}
\end{equation}
and hence,
\begin{equation}\nonumber
\displaystyle{\sum_{n=1}^{\infty}}n^{r+1}x^n\left(\frac{1-x}{x}\right) = \frac{P_{_r}(x)}{(1-x)^{r+1}}
\end{equation}
and thus proving that if eq.(\ref{genind}) is true for $k = r$, then it is also true for $k = r+1$. Thus by induction, we get eq.(\ref{genind}) is true for all $k$.
From eq.(\ref{genind}), we get
\begin{eqnarray}
E(Y^k) &=& \displaystyle{\sum_{n=1}^{\infty}}n^k \rho^n \left(\frac{1-\rho}{\rho}\right)\nonumber\\
       &=& \frac{P_{_{k-1}}(\rho)}{(1-\rho)^{k}}\nonumber
\end{eqnarray}
and hence
\begin{equation}
\label{ndiff8}
\displaystyle{\lim_{N \rightarrow \infty}}E(\underline{X}^k) = E(Y^k) = \displaystyle{\frac{P_{_{k-1}}(\rho)}{(1-\rho)^k}}
\end{equation}
eq.(\ref{ndiff8}) is a rational function of $\rho$ and hence a rational function of $R$. In particular, for the classical Fibonacci sequence, we have
\begin{equation}\nonumber
\displaystyle{\lim_{N \rightarrow \infty}} E(\underline{X}^k) = E(Y^k) = \frac{P_{_{k-1}}\left(\displaystyle{\frac{1}{\varphi}}\right)}{\left(1-\displaystyle{\frac{1}{\varphi}}\right)^k}
\end{equation}
Since $\displaystyle{\frac{1}{\varphi}} = \varphi-1$, we get
\begin{eqnarray}
\displaystyle{\lim_{N \rightarrow \infty}} E(\underline{X}^k) = E(Y^k) &=& \left(\frac{P_{_{k-1}}(\varphi-1)}{(\varphi-1)^k}\right)\varphi^k\nonumber\\
&=& \varphi^{2k} P_{_{k-1}}(\varphi-1)\nonumber\\
&=& Q(\varphi), \hspace{0.2in}\textup{a polynomial in}\hspace{0.1in} \varphi
\end{eqnarray}
But, since $\varphi^2 = \varphi + 1$, all powers $\varphi^j$ for $j\ge2$ can be expressed as linear polynomials in $\varphi$ and hence we get any polynomial in $\varphi$ can be rewritten as a linear polynomial in $\varphi$. Thus,
\begin{eqnarray}
 \displaystyle{\lim_{N \rightarrow \infty}} E(\underline{X}^k) &=& E(Y)\\
         &=& A_k\varphi + B_k
\end{eqnarray}
\end{proof}

We now provide a detailed analysis about the location of the maximum in the sequence resulting from self linear convolution of linear recurrence relations with positive integer sequences.

\section{Linear Convolution of linear recurrence relations with themselves}
In the previous sections, we proved certain optimal properties of the classical Fibonacci sequence among the general $k$-th order recurrence relations with positive integer coefficients, with respect to the discrete probability distributions that they induced. In this section we bring out another fascinating property of the classical Fibonacci sequence. Let $x_{_n}$ be a finite length integer sequence, of length $N$, generated by a $k$-th order recurrence relation with positive integer coefficients. The location of maxima in the sequence, $y_{_n}$ resulting from the linear convolution of $x_{_n}$ with itself, can be only one of the following:
\begin{itemize}
\item 2$N$-1, when the largest positive real root of the ``\textit{characteristic equation}'' is less than 2 
\item 2$N$-1, when the largest positive real root of the ``\textit{characteristic equation}'' is equal to  2 and the length even.
\item $2N$ when the largest positive real root is equal to 2 and the sequence length is odd.
\item $2N$ when the largest positive real root is greater than 2.
\end{itemize}
We prove the above in the following sections. 

\section{Some Notations and Definitions}
As mentioned earlier, while investigating the use of integer sequences for window functions in digital filtering \cite{arul}, we observed that most of the integer integer sequences generated by $k$-th order recurrence relations, with positive integer coefficients, are non-decreasing. However, conventional classical window functions are symmetric, bounded and smooth. Though, integer sequences, mentioned as above, are bounded for finite length, they are not symmetric. Hence to be used as window functions, integer sequences are symmetrized. However, this does not guarantee a smooth profile for most of the integer sequences. In order to get a smooth profile, we do one of the following:
\begin{itemize}
\item Perform linear convolution after symmetrizing the integer sequence about, say, $N/2$.
\item Perform linear convolution of the non-decreasing integer sequence and then symmetrically extend the resulting sequence about its maximum.
\end{itemize}
While in the first case, for any finite length, $N$, one can symmetrize at half the length, the second case requires the location of maximum in the sequence resulting from the linear convolution. In the sections to follow, we obtain the location of maxima in the sequence resulting from the linear convolution of $x_{_n}$ with itself. In particular, we look at the family of sequences generated using $k$-th order linear recurrence relation. First, we look at the second order case and then generalize the result to higher order recurrence relations.

\subsection{Notations}
Let
\begin{eqnarray}
 x &=& \{ x_{_n} \}_{n=1}^{^{N}}\\
 y &=& \{y_{_n}\}_{n=2}^{^{2N}} \\
\label{GenSecOrd} x_{_n} &=& ax_{_{n-1}} + bx_{_{n-2}}, \hspace{0.2in} \forall a,b \in \mathbb N 
\end{eqnarray}
The ``$characteristic$'' equation of the above is given by
\begin{equation}
x^2 - ax - b = 0
\end{equation}
The largest positive real root of the above equation is as below:
\begin{equation}
\label{LargestPositiveRealRoot}R = a + \displaystyle{\frac{\sqrt{a^2+4b}}{2}}
\end{equation}
 
\subsection{Definitions}
The linear convolution of two sequences is given by
\begin{eqnarray}
\label{FromntoN}y_{_n} = \displaystyle{\sum_{k=1}^{n-1}x_{_k} x_{_{n-k}}} \hspace{0.2in} 2 \leq n\leq N\\
\label{Fromnto2N}y_{_{N+j}} = \displaystyle{\sum_{k=j}^{N}x_{_k} x_{_{N+j-k}}} \hspace{0.2in} 1 \leq j \leq N
\end{eqnarray}

For $2\leq n \leq N$, we have,
\begin{eqnarray}
y_{_n} - y_{_{n-1}} &=& \displaystyle{\sum_{k=1}^{n-1}}x_{_k}x_{_{n-k}} -  \displaystyle{\sum_{k=1}^{n-2}}x_{_k}x_{_{n-1-k}}\\
&=& \displaystyle{\sum_{k=1}^{n-2}}x_{_k} (x_{_{n-k}}-x_{_{n-1-k}} )+ x_{_{n-1}}x_{_1}\\
& > & 0 \nonumber\\
\label{FornlessthanN}\Rightarrow y_{_n} &>& y_{_{n-1}} \hspace{0.2in} 1\leq n \leq N
\end{eqnarray}

From eq.(\ref{FornlessthanN}), we observe that the convolution sum is an increasing sequence for $1\leq n\leq N$.\\
Now we look at the eq.\ref{Fromnto2N}, for $j\geq 1$. 
\begin{equation}
\label{yNDifference}y_{_{N+j}} - y_{_{N+j-1}} = \displaystyle{\sum_{k=j}^{N}}x_{_{k}}x_{_{N+j-k}} -  \displaystyle{\sum_{k=j-1}^{N}}x_{_k}x_{_{N+j-1-k}}
\end{equation}
At this juncture, we revisit the general second order recurrence relation with positive integer coefficients, given in eq.(\ref{GenSecOrd}). We evaluate eq.(\ref{GenSecOrd}) for different possible values of $a$ and $b$ and hence compute the difference as in eq.(\ref{yNDifference}).
\begin{enumerate}
\item \label{a1} $a = 1$
	\begin{enumerate}
	\item \label{b1} $b = 1$
	\item \label{b2} $b = 2$
	\item \label{bg2}$b > 2$
	\end{enumerate}
\item \label{ag2} $a >1$
\end{enumerate}

\section{\underline{$a=1, b =1$}:}
We look at case(\ref{b1}) that corresponds to $a = b = 1$. Eq.(\ref{GenSecOrd}), in this case, generates the elements of classical Fibonacci sequence. Eq.(\ref{GenSecOrd}) can be rewritten as
\begin{equation}
x_{_n} = x_{_{n-1}} + x_{_{n-2}}
\end{equation}
$R$, the largest positive real root in the case of classical Fibonacci sequence corresponds to the golden ratio, $\varphi = \displaystyle{\left(\frac{1+\sqrt5}{2}\right)}$,
With the above conditions, eq.\ref{yNDifference} gives
\begin{eqnarray}
\label{FirstDiff}y_{_{N+j}} - y_{_{N+j-1}} &=& \displaystyle{\sum_{k=j}^{N}}x_{_{k}}x_{_{N+j-k}} -  \displaystyle{\sum_{k=j-1}^{N}}x_{_k}x_{_{N+j-1-k}}\\
&=& \displaystyle{\sum_{k=j}^{N}}x_{_{k}}(x_{_{N+j-1-k}}+x_{_{N+j-2-k}}) -  \displaystyle{\sum_{k=j-1}^{N}}x_{_k}x_{_{N+j-1-k}}\nonumber\\
&=& \displaystyle{\sum_{k=j}^{N}}x_{_{k}}x_{_{N+j-2-k}} - x_{_{j-1}}x_{_N}\nonumber
\end{eqnarray}
Suppose, $j \leq N-1$,
\begin{eqnarray}
y_{_{N+j}} - y_{_{N+j-1}}&\geq& x_{_{N-1}}x_{_{j-1}} + x_{_N}x_{_{j-2}} - x_{_{j-1}}x_{_N}\nonumber\\
&=& x_{_{N-1}}x_{_{j-1}} + x_{_N}(x_{_{j-1}}-x_{_{j-3}})- x_{_{j-1}}x_{_N}\nonumber\\
&=& x_{_{N-1}}x_{_{j-1}} - x_{_{j-3}}x_{_N}\nonumber\\
&=& x_{_{N-1}}(x_{_{j-2}}+x_{_{j-3}}) - x_{_{j-3}}(x_{_{N-1}}+x_{_{N-2}})\nonumber\\
\textup{Rearranging the terms, we get}\nonumber\\
&=& x_{_{N-1}}(x_{_{j-2}} - x_{_{j-3}}) + x_{_{j-3}}(x_{_{N-1}} - x_{_{N-2}})\nonumber
\end{eqnarray}
We find that the term on the right hand side, in the above equation is greater than 0. Hence
\begin{equation}
\label{jLessThan2Nminus1}y_{_{N+j}} - y_{_{N+j-1}} > 0
\end{equation}
From eq.(\ref{jLessThan2Nminus1}), we find the sequence $y_{_n}$ to be increasing for $1 \leq j \leq N-1$. The only pair of elements of $y_{_n}$ that needs to be compared are $y_{_{2N}}$ and $y_{_{2N-1}}$. From eq.(\ref{FirstDiff}), we get

For\hspace{0.1in} $N-1 \leq j \leq N$,\hspace{0.1in} the above corresponds to
\begin{eqnarray}
y_{_{2N}} - y_{_{2N-1}} &=& x_{_N}^2 - 2x_{_N}x_{_{N-1}}\\
&=& x_{_N}(x_{_N} - 2x_{_{N-1}})\nonumber\\
&=& x_{_N}(x_{_{N-2}}-x_{_{N-1}})\hspace{0.2in}\textup{since}\hspace{0.1in} x_{_N} = x_{_{N-1}} + x_{_{N-2}}\nonumber
\end{eqnarray} 
We know that \hspace{0.1in} $x_{{n}} < x_{_{n+1}}$, \hspace{0.1in} and hence, from the above equation, we have
\begin{equation}
\label{y2NDifference}\Rightarrow y_{_{2N}} - y_{_{2N-1}}< 0 
\end{equation}
From eq.(\ref{y2NDifference}), we find that $y_{_{2N}} < y_{_{2N-1}}$. Thus, in the sequence resulting from the linear convolution of finite length classical Fibonacci sequence with itself, the maximum is located at $2N-1$, where $N$ is the length of the sequence $x_{_n}$.

\section{\underline{$a=1, b = 2$}:}
\label{aequals1bequals2}
Now, we explore what happens when $a=1, b=2$ in eq.(\ref{GenSecOrd}), as mentioned in case(\ref{b2}). The corresponding second order recurrence relation is as below:
\begin{equation}
\label{a1b2Recurrence}x_{_n} =  x_{_{n-1}} + 2x_{{n-2}}
\end{equation}
In this case, the ``$characteristic$'' equation is $x^2 - x - 2 = 0$ and the largest positive root, $R$, is 2. To obtain the location of the maximum in the linear convolution of eq.\ref{a1b2Recurrence} with itself, we adopt the same procedure as classical Fibonacci sequence.
\begin{eqnarray}
y_{_{N+j}} - y_{_{N+j-1}} &=& \displaystyle{\sum_{k=j}^{N}}x_{_{k}}x_{_{N+j-k}} -  \displaystyle{\sum_{k=j-1}^{N}}x_{_k}x_{_{N+j-1-k}}\nonumber\\
&=& \displaystyle{\sum_{k=j}^{N}}x_{_{k}}(x_{_{N+j-1-k}}+2x_{_{N+j-2-k}}) -  \displaystyle{\sum_{k=j-1}^{N}}x_{_k}x_{_{N+j-1-k}}\nonumber\\
&=& \displaystyle{\sum_{k=j}^{N}}2x_{_{k}}x_{_{N+j-2-k}} \hspace{0.1in}- x_{_{j-1}}x_{_N}\\
\end{eqnarray}
Suppose, $j \leq N-1$,
\begin{eqnarray}
y_{_{N+j}} - y_{_{N+j-1}}&\geq& 2x_{_{N-1}}x_{_{j-1}} + x_{_N}x_{_{j-2}} - x_{_{j-1}}x_{_N}\nonumber\\
&\geq& 2x_{_{N-1}}x_{_{j-1}} + x_{_N}(x_{_{j-1}}-2x_{_{j-3}})- x_{_{j-1}}x_{_N}\nonumber\\
&=& 2x_{_{N-1}}x_{_{j-1}} - 2x_{_{j-3}}x_{_N}\nonumber\\
&=& 2x_{_{N-1}}(x_{_{j-2}}+2x_{_{j-3}}) - 2x_{_{j-3}}(x_{_{N-1}}+2x_{_{N-2}})\nonumber
\end{eqnarray}
Rearranging the terms, we get,
\begin{eqnarray}
&=& 2x_{_{N-1}}(x_{_{j-2}} - 2x_{_{j-3}}) + 4x_{_{j-3}}(x_{_{N-1}} - x_{_{N-2}})\nonumber\\
\label{a1b2jLessThan2Nminus1}y_{_{N+j}} - y_{_{N+j-1}}&>& 0
\end{eqnarray}
From eq.(\ref{a1b2jLessThan2Nminus1}), we find the sequence $y_{_n}$ to be an increasing one for $1 \leq j \leq N-1$. The only elements of $y_{_n}$ that need to be compared are $y_{_{2N}}$ and $y_{_{2N-1}}$.

For\hspace{0.1in} $N-1 \leq j \leq N$,\hspace{0.1in} the above corresponds to
\begin{eqnarray}
y_{_{2N}} - y_{_{2N-1}} &=& x_{_N}^2 - 2x_{_N}x_{_{N-1}}\\
\label{a1b2xnxnminus1}&=& x_{_N}(x_{_N} - 2x_{_{N-1}})
\end{eqnarray} 

\begin{itemize}
\item When $N$ is even, $x_{_N} < 2x_{_{N-1}}$ and hence $y_{_{2N-1}}$ is the maximum.
\item When $N$ is odd, $x_{_N} > 2x_{_{N-1}}$. Therefore, $y_{_{2N}}$ is maximum.
\end{itemize}
 This can be easily proved by mathematical induction.
\begin{proof}
Let $k$ = 2. With standard initial conditions 1, 1 one can observe that $f_2 < \frac{1}{2}f_3$ and $f_2 < 2f_1$.\\
Let $k$ = 4. The sequence elements with initial conditions 1, 1 turn out to be 1, 1, 3, 5. We observe that $f_4 < \frac{1}{2}f_5$ and $f_4 < 2f_3$. \\
For any $k$, we prove the following by induction. 
\begin{equation}
\label{induction}
f_{_{2k}} < \min\{2f_{_{2k-1}}, \frac{1}{2}f_{_{2k+1}}\}
\end{equation}
Let eq.(\ref{induction}) be true for $k$-1. So, we have
\begin{equation}
\label{inductionkminus1}
f_{_{2k-2}} < \min\{2f_{_{2k-3}}, \frac{1}{2}f_{_{2k-1}}\}
\end{equation}
We now prove that eq.(\ref{induction}) is true for $k$ i.e.,
\begin{equation}\nonumber
f_{_{2k}} < \min\{2f_{_{2k-1}}, \frac{1}{2}f_{_{2k+1}}\}
\end{equation}
\begin{eqnarray}
f_{_{2k}} &=& f_{_{2k-1}} + 2f_{_{2k-2}}\\
\textup{We~know~that~from~eq.(\ref{inductionkminus1})},f_{_{2k-2}} < \frac{1}{2}f_{_{2k-1}} \nonumber\\
\textup{Hence} f_{_{2k}} &<& f_{_{2k-1}} + f_{_{2k-1}}\nonumber\\
\Rightarrow \label{f2klessthan2f2kminus1} f_{_{2k}} &<& 2f_{_{2k-1}}
\end{eqnarray}
Similarly,
\begin{eqnarray}
f_{_{2k+1}} &=& f_{_{2k}} + 2f_{_{2k-1}}\\
\textup{From~eq.(\ref{f2klessthan2f2kminus1})},f_{_{2k}} < 2f_{_{2k-1}} \nonumber\\
\textup{Hence} f_{_{2k}} &<& f_{_{2k+1}} - f_{_{2k}} \nonumber\\
\Rightarrow 2f_{_{2k}} &<& f_{_{2k+1}} \nonumber\\
\Rightarrow f_{_{2k}} &<& \frac{1}{2}f_{_{2k+1}} 
\end{eqnarray}
Thus we find that, $f_{_{2k}} < \min\{2f_{_{2k-1}}, \frac{1}{2}f_{_{2k+1}}\}$. From this one can see that $y_{_n}$ has its maximum either at 2$N$-1 or 2$N$ depending on whether the sequence length is odd or even.
\end{proof}

\section{\underline{$a=1, b > 2$}:}
With $a=1, b > 2$ in eq.(\ref{GenSecOrd}), we look at case(\ref{bg2}). The corresponding second order recurrence relation is as below:
\begin{equation}
\label{a1bg2Recurrence}x_{_n} =  x_{_{n-1}} + bx_{{n-2}}
\end{equation}
In this case, the ``$characteristic$'' equation is $x^2 - x - b = 0$ and the largest positive root, $R$, greater than 2. We now proceed to obtain the location of the maximum in the linear convolution of eq.(\ref{a1bg2Recurrence}) with itself.
\begin{eqnarray}
y_{_{N+j}} - y_{_{N+j-1}} &=& \displaystyle{\sum_{k=j}^{N}}x_{_{k}}x_{_{N+j-k}} -  \displaystyle{\sum_{k=j-1}^{N}}x_{_k}x_{_{N+j-1-k}}\nonumber\\
&=& \displaystyle{\sum_{k=j}^{N}}x_{_{k}}(x_{_{N+j-1-k}}+b x_{_{N+j-2-k}}) -  \displaystyle{\sum_{k=j-1}^{N}}x_{_k}x_{_{N+j-1-k}}\nonumber\\
&=& \displaystyle{\sum_{k=j}^{N}}bx_{_{k}}x_{_{N+j-2-k}} \hspace{0.1in}- x_{_{j-1}}x_{_N}\\~\nonumber
\end{eqnarray}
Suppose, $j \leq N-1$,
\begin{eqnarray}
y_{_{N+j}} - y_{_{N+j-1}}&\geq& bx_{_{N-1}}x_{_{j-1}} + b x_{_N}x_{_{j-2}} - x_{_{j-1}}x_{_N}\nonumber\\
&\geq& bx_{_{N-1}}x_{_{j-1}} + b x_{_N}(x_{_{j-1}}- b x_{_{j-3}})- x_{_{j-1}}x_{_N}\nonumber\\
&=& bx_{_{N-1}}x_{_{j-1}} - bx_{_{j-3}}x_{_N}\nonumber\\
&=& bx_{_{N-1}}(x_{_{j-2}} + bx_{_{j-3}}) - b x_{_{j-3}}(x_{_{N-1}} + bx_{_{N-2}})\nonumber
\end{eqnarray}
Rearranging the terms, we get,
\begin{eqnarray}
&=& bx_{_{N-1}}(x_{_{j-2}} - x_{_{j-3}}) + b^2x_{_{j-3}}(x_{_{N-1}} - x_{_{N-2}})\\
\label{a1bg2jLessThan2Nminus1}y_{_{N+j}} - y_{_{N+j-1}}&>& 0
\end{eqnarray}
From eq.(\ref{a1bg2jLessThan2Nminus1}), we find the sequence $y_{_n}$ to be an increasing one for $1 \leq j \leq N-1$. The only elements of $y_{_n}$ that need to be compared are $y_{_{2N}}$ and $y_{_{2N-1}}$.

For\hspace{0.1in} $N-1 \leq j \leq N$,\hspace{0.1in} the above corresponds to
\begin{eqnarray}
y_{_{2N}} - y_{_{2N-1}} &=& x_{_N}^2 - 2x_{_N}x_{_{N-1}}\\
&=& x_{_N}(x_{_N} - 2x_{_{N-1}})\\
\label{a1bg2testequation}&=& x_{_N}(x_{_{N-1}}+ b x_{_{N-2}}- 2x_{_{N-1}})
\end{eqnarray} 
Let us consider the sequence generated by eq.(\ref{a1bg2Recurrence}), with $b=3$. The elements of the sequence, with initial values being 1, 1, are
\begin{center}
1, 1, 4, 7, 19, 40, 97, $\ldots$
\end{center}
From the above, we find that $b x_{_{n-2}}$ is always greater than $x_{_{n-1}}$. Thus eq.\ref{a1bg2testequation} is always positive. Hence, we find that $y_{_{2N}} > y_{_{2N-1}}$. Thus, the maximum is located at 2$N$ and not at 2$N-1$. A formal proof by induction can be obtained along the same lines of mathematical induction, as found in section. \ref{aequals1bequals2}.

\section*{\underline{{$a>1$}}:}
We now present the case when $a > 1$ and $b \in \mathbb N$, in eq.(\ref{GenSecOrd}).
\begin{eqnarray}
y_{_{N+j}} - y_{_{N+j-1}} &=& \displaystyle{\sum_{k=j}^{N}}x_{_{k}}x_{_{N+j-k}} -  \displaystyle{\sum_{k=j-1}^{N}}x_{_k}x_{_{N+j-1-k}}\nonumber\\
&=& \displaystyle{\sum_{k=j}^{N}}x_{_{k}}(a x_{_{N+j-1-k}}+b x_{_{N+j-2-k}}) -  \displaystyle{\sum_{k=j-1}^{N}}x_{_k}x_{_{N+j-1-k}}\nonumber\\
&=& \displaystyle{\sum_{k=j}^{N}}(a-1)x_{_{k}}x_{_{N+j-1-k}} + \displaystyle{\sum_{k=j}^{N}}bx_{_{k}}x_{_{N+j-2-k}} \hspace{0.1in}- x_{_{j-1}}x_{_N}\\
\end{eqnarray}
When $a >1$, the term $x_{_{j-1}}x_{_N}$ gets accounted for by $\displaystyle{\sum_{k=j}^{N}}(a-1)x_{_{k}}x_{_{N+j-1-k}}$, with $k=N$. The remaining terms are all positive. Hence $y_{_n}$ increases till $2N-1$. Now, we compare $y_{_{2N-1}}$ and $y_{_{2N}}$.
For\hspace{0.1in} $N-1 \leq j \leq N$,\hspace{0.1in} the above corresponds to
\begin{eqnarray}
y_{_{2N}} - y_{_{2N-1}} &=& x_{_N}^2 - 2x_{_N}x_{_{N-1}}\\
&=& x_{_N}(x_{_N} - 2x_{_{N-1}})\\
&=& x_{_N}(a x_{_{N-1}}+ b x_{_{N-2}}- 2x_{_{N-1}})\\
\label{ag1btestequation}&=& x_{_N}((a-2) x_{_{N-1}}+ b x_{_{N-2}})
\end{eqnarray} 
With $a > 1$ and $b \in \mathbb N$, eq.(\ref{ag1btestequation}) is always positive. Hence $y(n)$ has its maximum at $y_{_{2N}}$.

\section{Self linear convolution of higher order linear recurrence relations with integer coefficients}
In the previous section, we found that the sequence resulting from the self linear convolution of second order linear recurrence (with positive integer coefficients) has its maximum located at either 2$N$-1 or 2$N$ is determined by the
\begin{itemize} 
\item coefficients
\item and hence the largest positive real root. 
\end{itemize}
It was proved that, when the largest positive real root was 2, the location of maximum is either 2$N$-1 or 2$N$, based on whether $N$, the length of $x_{_n}$, is even or odd.
 It is also worth noting that this is true even in the case of self convolution of, $x_{_n} = a x_{_{n-1}}$, first order linear recurrence relation (with positive integer coefficients). In this case, if $a$ = 2, the largest positive real root is 2.  $y_{_n}$ has its maximum at either 2$N$-1 or 2$N$ depending on whether the length $N$ is even or odd.  However, when $a$ = 1, it can be easily proved that, the maximum of $y_{_n}$ is located at $N$. \\
Now, we look at the higher order recurrence relations with positive integer coefficients. Let 
\begin{equation}
 x_{_n} = a_{_{n-1}}x_{_{n-1}} + \ldots + a_{_k}x_{_{n-k}} \hspace{0.3in} \forall a_{_i} \in \mathbb{N}
\end{equation}
be a higher order recurrence relation. On similar lines, as self linear convolution of the second order linear recurrence relation, one can establish that even in the case of higher order recurrence relations, there exist three cases, namely,

\begin{itemize}
\item all the $a_{_i}$'s are 1. $y_{_n}$ is maximum at $n = 2N-1$
\item $a_{_i}$'s are 1 and $a_{_0}$ = 2. $y_{_n}$ is maximum at $n = 2N-1$ or $n = 2N$, depending on whether the length is even or odd
\item $a_{_i}$'s are greater than 1. $y_{_n}$ is maximum at $n = 2N$
\end{itemize}

\section{Conclusion}
In this work, on the one hand, we have explored and proved certain optimal probabilistic limit properties of Fibonacci sequence. We have proved that, among the distributions induced by $k^{th}$ order linear recurrence equations with positive integer coefficients, the limits of the mean and the variance are maximum for those generated by second order linear recurrence relation, corresponding to the Fibonacci sequence. We also found that the ratio of the variance to the mean, of the distribution induced by the Fibonacci sequence turns out to be golden ratio, $\varphi$. Moreover, for all such distributions, we have,
\begin{itemize}
\item $\underline{X}$ converges to $Y$ in distribution
\item all the moments of $\underline{X}$ converge to the corresponding moments of $Y$
\item in the classical Fibonacci case, all moments of $\underline{X}$ converges to linear function of the golden ratio, $\varphi$.
\end{itemize}
Besides proving certain optimal probabilistic limit properties, we also analyzed the second order linear recurrence relation with positive integer coefficients, from the point of view of self linear convolution. We established that the sequence resulting from the self linear convolution of any higher order linear recurrence relation with integer coefficients, can have its maximum either at 2$N$ or 2$N$-1. A more interesting observation in this regard is that when the largest positive real root of such recurrence relations is 2, the location of maximum depends on whether the sequence length of $x_{_n}$ is odd or even. 

\appendix
\section{Some Preliminaries}
\label{Prelims}
We list some of the identities used in this paper in the following sections.
\begin{equation}
\label{SumxN} G(x) = \displaystyle{\sum_{n=1}^{N}}x^n = x+x^2+\ldots+x^N = \displaystyle{\frac{x^{N+1}-x}{x-1}}
\end{equation}
Let $G'(x)$ be the derivative of $G(x)$. We have,
\small
\begin{eqnarray}
G'(x) = 1+2x+\ldots+Nx^{N-1} &=& \displaystyle{\frac{(x-1)\{(N+1)x^{N}-1\}- (x^{N+1}-x)}{(x-1)^2)}}\nonumber\\
     &=& \displaystyle{\frac{(N+1)x^{N+1}-x-(N+1)x^{N}+1-x^{N+1}+x}{(x-1)^2}}\nonumber\\
     &=& \displaystyle{\frac{Nx^{N+1}-(N+1)x^N + 1}{(x-1)^2}}\nonumber
\end{eqnarray}
\normalsize
Hence, we get,
\begin{eqnarray}
xG'(x) = \displaystyle{\sum_{n=1}^{N}}nx^n &=& x+2x^2+\ldots+ Nx^N \nonumber\\
\label{SumNxN}&=&\displaystyle{\frac{Nx^{N+2}-(N+1)x^{N+1}+x}{(x-1)^2}}
\end{eqnarray}
From this we get,
\begin{equation}
(xG'(x))' = 1+2^2x+\ldots+N^2x^{N-1}
\end{equation}
which by eq.(\ref{SumNxN}) is 
\scriptsize
\begin{eqnarray}
\nonumber(xG'(x))' &=& \displaystyle{\frac{(x-1)^2\{N(N+2)x^{N+1}-(N+1)^2x^N+1\}-2(x-1)\{Nx^{N+2}-(N+1)x^{N+1}+x\}}{(x-1)^4}}\\
\nonumber&=& \displaystyle{\frac{(x-1)\{N(N+2)x^{N+1}-(N+1)^2x^N+1\}-2\{Nx^{N+2}-(N+1)x^{N+1}+x\}}{(x-1)^3}}\\
\nonumber&=& \displaystyle{\frac{N^2\{x^{N+2}-2x^{N+1}+x^N\}+ N\{-2x^{N+1}+2x^N\}+\{x^{N+1}+x^N+x+1\}}{(x-1)^3}}
\end{eqnarray}
\normalsize
and hence,
\scriptsize
\begin{eqnarray}
\nonumber x(xG'(x))' &=& \displaystyle{\sum_{n=1}^{N}}n^2x^n \\
\label{SumN2xN}&=& \displaystyle{\frac{N^2\{x^{N+3}-2x^{N+2}+x^{N+1}\}+ N\{-2x^{N+2}+2x^{N+1}\}+\{x^{N+2}+x^{N+1}+x^2+x\}}{(x-1)^3}}
\end{eqnarray}
\normalsize

\section{Useful Finite Sums}
\label{FiniteSums}
In this section, we obtain the expressions for some of the finite sums relevant to the work presented in this paper. 
\subsection{$\displaystyle{\sum_{n=1}^Nf[n]}$}
Let $\bar{S}$ be the sum of the sequence $f[n]$ defined as
\begin{equation}
\label{fsum}
\bar{S} =\displaystyle{\sum_{n=1}^N f[n]} = \displaystyle{\sum_{n=1}^N \frac{1}{\sqrt{a^2 + 4b}}(R^n - R_S^n)}
\end{equation}
From eq.(\ref{SumxN}), we get
\begin{equation}
\label{sumGR}
G(R) = R + R^2 + R^3 + \ldots + R^N =\displaystyle{ \sum_{n=1}^N R^n = \frac{R^{N+1}- R}{R-1}}
\end{equation}
\begin{equation}
\label{sumGR1}
G(R_S) = R_S + R_S^2 + R_S^3 + \ldots + R_S^N = \displaystyle{\sum_{n=1}^N R_S^n = \frac{R_S^{N+1}- R_S}{R_S-1}}
\end{equation}
Substituting eq.(\ref{sumGR}) and eq.(\ref{sumGR1}) eq.(\ref{fsum}), we get,
\begin{eqnarray}
\label{GRG1R}
\bar{S} =& \displaystyle{\frac{1}{\sqrt{a^2 + 4b}}}\displaystyle{\sum_{n=1}^N{(G(R) -  G(R_S))}}\nonumber\\
  \label{ratio} = & \displaystyle{\frac{1}{\sqrt{a^2 + 4b}}} \left(\displaystyle{\frac{R^{N+1}- R}{R-1}} - \displaystyle{\frac{R_S^{N+1}- R_S}{R_S-1}}\right)
\end{eqnarray}
We have,
\begin{eqnarray}
\nonumber R^2 &=& \displaystyle{\frac{2a^2+4b+2a\sqrt{a^2+4b}}{4}}\\ 
\nonumber \Rightarrow R^2 &>& b\\
\Rightarrow |\frac{R_S}{R}| &=& \frac{b}{R^2} < 1  
\end{eqnarray}
Hence, taking $R^{N+1}$ common from eq.(\ref{ratio}), we see that,
\begin{equation}
\label{SumIncrease}
 \bar{S} \sim \displaystyle{\frac{1}{\sqrt{a^2 + 4b}}}\left( \frac{R^{N+1}}{R-1}\right)
\end{equation}

\subsection{$\displaystyle{\sum_{n=1}^{N}n f[n]}$}
Let $\bar{S}_1$ be defined as 

\begin{eqnarray}
\label{nfn}\bar{S}_1 &=& \displaystyle{\sum_{n=1}^N nf[n]} \nonumber\\
&=& \frac{1}{\sqrt{a^2 + 4b}} \displaystyle{\sum_{n=1}^N{(n R^n -  n(R_S)^n)}}\nonumber
\end{eqnarray}
By eq.(\ref{SumNxN}), we get
\begin{equation}
\displaystyle{\sum_{n=1}^N }nR^n = \frac{NR^{N+2}-(N+1)R^{N+1}+R}{(R-1)^2}
\end{equation}
We can get a similar expression involving the other root $R_S$. Using the fact that the absolute value of the ratio of $R_S$ and $R$ is less than 1, we obtain 
\begin{equation}
\label{nfn1}
\bar{S}_1 \sim \frac{1}{\sqrt{a^2 + 4b}}\left(\frac{NR^{N+2}-(N+1)R^{N+1}+R}{(R-1)^2}\right)
\end{equation}

\subsection{$\displaystyle{\sum_{n=1}^{N}n^2 f[n]}$}
Let $\bar{S}_2$ be defined as
 \begin{equation}
 \label{n2fn}
 \bar{S}_2 = \displaystyle{\sum_{n=1}^{N}n^2 f[n]}
 \end{equation}
 Using eq.(\ref{SumN2xN}) we get 
 \begin{equation}
 \label{s2}
 \bar{S}_2 \sim \frac{N^2\{R^{N+3}-2R^{N+2}+R^{N+1}\}+N\{-2R^{N+2}+2R^{N+1}\}+ \{R^{N+2}+R^{N+1}+R^2+R\}}{(R-1)^3}
\end{equation}

\subsection{$\displaystyle{\sum_{n=1}^{N}f[N+1-n]}$}
Let $\underline{S}$ be defined as
\begin{eqnarray}
\label{sumd}
\underline{S} &=& \displaystyle{\sum_{n=1}^{N}f[N+1-n]}\\
\label{SumDecrease}\textup{Clearly,} \hspace{0.3in} \bar{S} = \displaystyle{\sum_{n=1}^{N}f[n]} &=& \displaystyle{\sum_{n=1}^{N}f[N+1-n]} = \underline{S}
\end{eqnarray}
Thus the sum $\underline{S}$ = $\bar{S}$ , where $\bar{S}$ is as in eq.(\ref{SumIncrease}).

\subsection{$\displaystyle{\sum_{n=1}^{N}nf[N+1-n]}$}
Let $\underline{S}_1$ denote $\displaystyle{\sum_{n=1}^{N}nf[N+1-n]}$
\begin{eqnarray}
\underline{S}_1 = \displaystyle{\sum_{n=1}^{N}nf[N+1-n]}\nonumber\\
= -\displaystyle{\sum_{n=1}^{N}(-n)f[N+1-n]}\nonumber\\
= - \displaystyle{\sum_{n=1}^{N}(N+1-n)f[N+1-n]} + (N+1)\displaystyle{\sum_{n=1}^{N}f[N+1-n]}\nonumber\\
= - \displaystyle{\sum_{k=1}^{N}kf[k] + (N+1)\bar{S}}\nonumber\\
\label{nfnd}\Rightarrow  \underline{S}_1 = (N+1)\bar{S} - \bar{S}_1
\end{eqnarray}

\subsection{$\displaystyle{\sum_{n=1}^{N}n^2f[N+1-n]}$}
Let $\underline{S}_2$ be defined as
\begin{eqnarray}
\underline{S}_2  = \displaystyle{\sum_{n=1}^{N}n^2f[N+1-n]}\nonumber\\
= \displaystyle{\sum_{n=1}^{N}(N+1-n)^2f[n]}\nonumber\\
= (N+1)^2 \displaystyle{\sum_{n=1}^Nf[n]} - 2(N+1)\displaystyle{\sum_{n=1}^N nf[n]} + \displaystyle{\sum_{n=1}^N n^2f[n]}\nonumber\\
\label{n2fnd}\Rightarrow \underline{S}_2 = (N+1)^2\bar{S} -2(N+1)\bar{S}_1 + \bar{S}_2
\end{eqnarray}
where $\bar{S}_1$ and $\bar{S}_2$ are as in eq.(\ref{nfn}) and in eq.(\ref{n2fn}) respectively.

\medskip

\noindent MSC2010: 11B37, 11B39, 60C05, 60C09, 60F99

\end{document}